\documentclass[12pt,leqno]{article}
\usepackage{amsfonts}
\usepackage{amsmath, amsthm, amsfonts, amssymb, color}
\usepackage{mathrsfs}
\usepackage{color}
\theoremstyle{definition}

\newcommand{\scr}[1]{\mathscr #1}
\definecolor{wco}{rgb}{0.5,0.2,0.3}

\numberwithin{equation}{section} \theoremstyle{remark}

\newcommand{\ua}{\uparrow}

\title{{\bf    Exponential Convergence in Entropy and Wasserstein  for McKean-Vlasov SDEs }\footnote{Supported in
 part by  NNSFC (11771326, 11831014, 11921001) and the DFG through the CRC 1283.} }
\author{
{\bf Panpan Ren$^{c)}$,    Feng-Yu Wang$^{a,b)}$  }\\
\footnotesize{$^{a)}$ Center for Applied Mathematics, Tianjin University, Tianjin 300072, China}\\
 \footnotesize{$^{b)}$ Department of Mathematics,
Swansea University, Bay Campus, SA1 8EN, United Kingdom}\\
\footnotesize{$^{c)}$ Department of Mathematics,
City University of HongKong, HongKong, China}\\
\footnotesize{ rppzoe@gmail.com,  wangfy@tju.edu.cn}}
\begin{document}
\allowdisplaybreaks
\def\R{\mathbb R}  \def\ff{\frac} \def\ss{\sqrt} \def\B{\mathbf
B}
\def\N{\mathbb N} \def\kk{\kappa} \def\m{{\bf m}}
\def\ee{\varepsilon}\def\ddd{D^*}
\def\dd{\delta} \def\DD{\Delta} \def\vv{\varepsilon} \def\rr{\rho}
\def\<{\langle} \def\>{\rangle}
  \def\nn{\nabla} \def\pp{\partial} \def\E{\mathbb E}
\def\d{\text{\rm{d}}} \def\bb{\beta} \def\aa{\alpha} \def\D{\scr D}
  \def\si{\sigma} \def\ess{\text{\rm{ess}}}\def\s{{\bf s}}
\def\beg{\begin} \def\beq{\begin{equation}}  \def\F{\scr F}
\def\Ric{\mathcal Ric} \def\Hess{\text{\rm{Hess}}}
\def\e{\text{\rm{e}}} \def\ua{\underline a} \def\OO{\Omega}  \def\oo{\omega}
 \def\tt{\tilde}\def\[{\lfloor} \def\]{\rfloor}
\def\cut{\text{\rm{cut}}} \def\P{\mathbb P} \def\ifn{I_n(f^{\bigotimes n})}
\def\C{\scr C}      \def\aaa{\mathbf{r}}     \def\r{r}
\def\gap{\text{\rm{gap}}} \def\prr{\pi_{{\bf m},\varrho}}  \def\r{\mathbf r}
\def\Z{\mathbb Z} \def\vrr{\varrho} \def\ll{\lambda}
\def\L{\scr L}\def\Tt{\tt} \def\TT{\tt}\def\II{\mathbb I}
\def\i{{\rm in}}\def\Sect{{\rm Sect}}  \def\H{\mathbb H}
\def\M{\mathbb M}\def\Q{\mathbb Q} \def\texto{\text{o}} \def\LL{\Lambda}
\def\Rank{{\rm Rank}} \def\B{\scr B} \def\i{{\rm i}} \def\HR{\hat{\R}^d}
\def\to{\rightarrow}\def\l{\ell}\def\iint{\int}\def\gg{\gamma}
\def\EE{\scr E} \def\W{\mathbb W}
\def\A{\scr A} \def\Lip{{\rm Lip}}\def\S{\mathbb S}
\def\BB{\scr B}\def\Ent{{\rm Ent}} \def\i{{\rm i}}\def\itparallel{{\it\parallel}}
\def\g{{\mathbf g}}\def\Sect{{\mathcal Sec}}\def\T{\mathcal T}\def\BB{{\bf B}}
\def\f{\mathbf f} \def\g{\mathbf g}\def\BL{{\bf L}}  \def\BG{{\mathbb G}}
\def\Bd{{D^E}} \def\BdP{D^E_\phi} \def\Bdd{{\bf \dd}} \def\Bs{{\bf s}} \def\GA{\scr A}
\def\Bg{{\bf g}}  \def\Bdd{\psi_B} \def\supp{{\rm supp}}\def\div{{\rm div}}
\def\ddiv{{\rm div}}\def\osc{{\bf osc}}\def\1{{\bf 1}}\def\BD{\mathbb D}\def\GG{\Gamma}
\maketitle

\begin{abstract}  The following type of exponential convergence  is proved for  (non-degenerate or degenerate) McKean-Vlasov  SDEs:
$$\W_2(\mu_t,\mu_\infty)^2 +\Ent(\mu_t|\mu_\infty)\le c \e^{-\ll t} \min\big\{\W_2(\mu_0, \mu_\infty)^2,\Ent(\mu_0|\mu_\infty)\big\},\ \ t\ge 1,$$
where $c,\ll>0$ are constants, $\mu_t$ is the distribution of the solution at time $t$, $\mu_\infty$ is the unique invariant probability measure, $\Ent$ is the relative entropy and $\W_2$ is the $L^2$-Wasserstein distance. In particular, this type of exponential convergence holds for  some (non-degenerate or degenerate) granular media type equations  generalizing those studied in \cite{CMV, GLW}  on the exponential convergence in  a mean field entropy. \end{abstract} \noindent

 \vskip 2cm

 \section{Introduction}

 The convergence in entropy for stochastic systems is an important topic in both probability theory and  mathematical physics, and has been well studied for Markov processes by using the log-Sobolev inequality, see for instance \cite{BGL} and references therein. However, the existing results derived in the literature do not apply to McKean-Vlasov SDEs due to the nonlinearity of the associated Fokker-Planck equations.
   In 2003, Carrillo, McCann and Villani \cite{CMV} proved  the exponential convergence in a mean field entropy  of the  following granular media equation for   probability density functions $(\rr_t)_{t\ge 0}$ on $\R^d$:
 \beq\label{E0} \pp_t \rr_t= \DD\rr_t  + {\rm div} \big\{\rr_t\nn (V + W*\rr_t)\big\},\end{equation}
 where the internal potential $V\in C^2( \R^d)$  satisfies   $\Hess_V\ge \ll  I_{d}$ for a constant $\ll>0$ and the $d\times d$-unit matrix $I_d$,  and the interaction potential $W\in C^2(\R^d)$ satisfies
 $W(-x)=W(x)$ and  $\Hess_W\ge -\dd I_{d}$
 for some constant $\dd\in [0, \ll/2)$.   Recall  that  we write $M\ge \ll I_d$ for a constant $\ll$ and a $d\times d$-matrix $M$, if $\<Mv,v\>\ge \ll |v|^2$ holds for any $v\in \R^d$.
 To introduce the mean field entropy, let  $\mu_V(\d x):= \ff{\e^{-V(x)}\d x}  {\int_{\R^d}\e^{-V(x)}\d x}$,   recall  the classical relative entropy
 $${\rm Ent}(\nu|\mu) := \beg{cases} \mu(\rr\log\rr),    &\text{if} \ \nu=\rr\mu,\\
  \infty, &\text{otherwise}\end{cases}$$ for $\mu,\nu\in \scr P$,   the space  of all probability measures on $\R^d$, and consider  the free energy functional
  $$E^{V,W}(\mu):= {\rm Ent}(\mu|\mu_V)+ \ff 1 2 \int_{\R^d\times\R^d} W(x-y) \mu(\d x)\mu(\d y),\ \ \mu\in \scr P,$$
  where we set $E^{V,W}(\mu)=\infty$ if either ${\rm Ent}(\mu|\mu_V)=\infty$ or the integral term is not well defined.  Then the
 associated mean field entropy ${\rm Ent}^{V,W}$
   is defined by
\beq\label{ETP}  {\rm Ent}^{V,W}(\mu):= E^{V,W}(\mu) -  \inf_{\nu\in\scr P}  E^{V,W}(\nu),\ \ \mu\in \scr P.\end{equation}
According to \cite{CMV},  for $V$ and $W$ satisfying the above mentioned conditions,   $E^{V,W}$ has a unique minimizer $\mu_\infty,$ and  $\mu_t(\d x):= \rr_t(x)\d x$  for probability density $\rr_t$ solving \eqref{E0} converges to $\mu_\infty$  exponentially    in the mean field  entropy:
$${\rm Ent}^{V,W}(\mu_t)\le \e^{-(\ll-2\dd)t} {\rm Ent}^{V,W}(\mu_0),\ \ t\ge 0.$$
Recently, this result was generalized  in \cite{GLW} by establishing the uniform log-Sobolev inequality for the associated mean field particle systems, such that   ${\rm Ent}^{V,W}(\mu_t)$ decays exponentially   for a class of non-convex $V\in C^2(\R^d)$ and $W\in C^2(\R^d\times \R^d)$, where  $W(x,y)=W(y,x)$ and  $\mu_t(\d x):=\rr_t(x)\d x$ for $\rr_t$ solving the nonlinear PDE
\beq\label{E00} \pp_t \rr_t= \DD\rr_t  + {\rm div} \big\{\rr_t\nn (V + W\circledast\rr_t)\big\},\end{equation}
where
\beq\label{AOO}  W\circledast\rr_t:=\int_{\R^d}W(\cdot,y) \rr_t(y)\d y.\end{equation}     In this case,   $\Ent^{V,W}$ is defined in \eqref{ETP}  for  the free energy functional
$$E^{V,W}(\mu):=  {\rm Ent}(\mu|\mu_V)+ \ff 1 2 \int_{\R^d\times\R^d} W(x,y) \mu(\d x)\mu(\d y),\ \ \mu\in \scr P.$$
To study \eqref{E00} using probability methods, we consider the
     following McKean-Vlasov SDE with initial distribution $\mu_0$:
\beq\label{E0'} \d X_t= \ss 2 \d B_t-\nn \big\{ V+ W\circledast\L_{X_t}\big\}(X_t)\d t,\end{equation}
where $B_t$ is the $d$-dimensional Brownian motion,  $\L_{X_t}$ is the distribution of $X_t$, and
\beq\label{A01} (W\circledast\mu)(x):=\int_{\R^d} W(x,y) \mu(\d y),\ \ x\in \R^d, \mu\in \scr P\end{equation}   provided the integral exists.
Let $\rr_t(x)=\ff{(\L_{X_t})(\d x)}{\d x},\ t\ge 0.$ By It\^o's formula and the integration by parts formula,   we have
\beg{align*}& \ff{\d}{\d t} \int_{\R^d} ( \rr_t f)(x)\d x =\ff{\d }{\d t} \E [f(X_t)] =\E\big[\big(\DD-\nn V- \nn \{W\circledast\rr_t\}\big)f(X_t)\big]\\
&= \int_{\R^d} \rr_t(x) \big\{\DD f -\<\nn V+ \nn \{W\circledast\rr_t\}, \nn f\>  \big\}(x)\d x\\
&= \int_{\R^d} f(x) \{\DD\rr_t + {\rm div}[\rr_t\nn V+ \rr_t \nn (W\circledast\rr_t)]\big\}(x)\d x,\ \ t\ge 0,\ f\in C_0^\infty(\R^d).\end{align*}
Therefore, $\rr_t$ solves \eqref{E00}. On the other hand, by this fact and the uniqueness of \eqref{E0} and \eqref{E0'}, if $\rr_t$ solves \eqref{E0} with $\mu_0(\d x):=\rr_0(x)\d x$,
then $\rr_t(x)\d x=\L_{X_t}(\d x)$ for $X_t$ solving \eqref{E0'} with $\L_{X_0}=\mu_0.$

To extend the study of \cite{CMV, GLW}, in this paper we   investigate the exponential convergence in entropy for  the following McKean-Vlasov SDE on $\R^d$:
\beq\label{E1} \d X_t= \si(X_t)\d W_t+ b(X_t,\L_{X_t})\d t,\end{equation}
where $W_t$ is the $m$-dimensional Brownian  motion on a complete filtration probability space $(\OO, \{\F_t\}_{t\ge 0}, \P)$,
$$\si: \R^d\to \R^{d}\otimes \R^m,\ \ b:\R^d\times \scr P_2\to \R^d$$ are measurable, and $\scr P_2$ is the class of probability measures on $\R^d$ with $\mu(|\cdot|^2)<\infty$.

Since   the  $``$mean field  entropy" associated with the SDE \eqref{E1} is not available, and   it is less explicit even exists  as  in \eqref{ETP}   for  the special model \eqref{E0'},   we intend to study the exponential convergence of $\L_{X_t}$ in the classical relative entropy $ {\rm Ent}$ and the Wasserstein distance $\W_2$.
Recall that for any $p\ge 1$,   the $L^p$-Wasserstein distance is defined by
$$\W_p(\mu_1,\mu_2):= \inf_{\pi\in \C(\mu_1,\mu_2)} \bigg(\int_{\R^d\times \R^d} |x-y|^p\pi(\d x,\d y)\bigg)^{\ff 1 p},\ \ \mu_1,\mu_2\in \scr P_p,$$
where $\C(\mu_1,\mu_2)$ is the set of all couplings of $\mu_1$ and $\mu_2$.

Unlike in \cite{CMV, GLW} where the mean field particle systems are used to estimate the mean field entropy, in this paper we use 
the log-Harnack inequality introduced in \cite{W10, RW10} and the Talagrand inequality developed in \cite{TAL, BGL, OV},  see Theorem \ref{T0} below. Thus, the key point of the present study is to establish these two types of inequalities for McKean-Vlasov SDEs.
Since the log-Harnack inequality is not yet available when $\si$  depends on the distribution,   in \eqref{E1} we only consider  distribution-free $\si$.   In particular, for  a class of   granular media type equations  generalizing the framework of \cite{CMV,GLW},  we prove
$$ \W_2(\mu_t,\mu_\infty)^2 +\Ent(\mu_t|\mu_\infty)\le c \e^{-\ll t} \min\big\{\W_2(\mu_0, \mu_\infty)^2,\Ent(\mu_0|\mu_\infty)\big\},\ \ t\ge 1$$
for $\mu_t(\d x):= \rr_t(x)\d x$ and some constants $c,\ll>0$,  see Theorem  \ref{C2.0}  below for details.

\

The remainder of the paper is organized as follows. In Section 2,  we state our main results for non-degenerate and degenerate models respectively, where  the first case includes the granular media type equations  \eqref{E00} or the corresponding Mckean-Vlasov SDE \eqref{E0'} as a special example, and the second case deals with  the McKean-Vlasov stochastic Hamiltonian system  referring to the degenerate granular media equation.
The main  results are   proved in Sections 3-5 respectively, where  Section 4  establishes  the log-Harnack inequality  for McKean-Vlasov stochastic Hamiltonian systems.

\section{Main results and examples }

We first present a criterion on the exponential convergence for McKean-Vlasov  SDEs by using the log-Harnack and Talagrand  inequalities, and prove  \eqref{WU} for  the granular media type equations \eqref{E01} below which generalizes  the framework of \cite{GLW}.  Then we state our results for solutions of SDE \eqref{E1} with non-degenerate and degenerate noises respectively.

\subsection{A criterion with  application to Granular media type equations}

In general, we consider the following McKean-Vlasov   SDE:
\beq\label{EP}  \d X_t= \si(X_t,\L_{X_t})\d W_t+ b(X_t,\L_{X_t})\d t,\end{equation}
where $W_t$ is the $m$-dimensional Brownian motion and
$$\si: \R^d\times \scr P_2\to \R^{d}\otimes \R^m,\ \ b:\R^d\times\scr P_2 \to \R^d$$ are measurable.  We assume that this SDE is strongly and weakly well-posed
for square integrable initial values. It is in particular the case if $b$ is continuous on $\R^d\times\scr P_2$ and there exists a constant $K>0$ such that
\beq\label{KK}\beg{split} &\<b(x,\mu)-b(y,\nu), x-y\>^+ +\|\si(x,\mu)-\si(y,\nu)\|^2\le K \big\{|x-y|^2 +\W_2(\mu,\nu)^2\},\\
& |b(0,\mu)|\le c\Big(1+\ss{\mu(|\cdot|^2)}\Big),\ \ x,y\in \R^d, \mu,\nu\in \scr P_2,\end{split}\end{equation}
see for instance \cite{W18}. See also \cite{HW20,Zhao20} and references therein for the well-posedness of McKean-Vlasov  SDEs with singular coefficients. For any  $\mu\in \scr P_2$, let $P_t^*\mu=\L_{X_t}$ for the solution $X_t$ with initial distribution $\L_{X_0}=\mu$.  Let
$$P_t f(\mu)= \E[f(X_t)]=\int_{\R^d}f\d P_t^*\mu,\ \ t\ge 0, f\in \B_b(\R^d).$$
We have the following equivalence on the exponential convergence of $P_t^*\mu$ in $\Ent$ and $\W_2$.

\beg{thm}\label{T0} Assume that $P_t^*$ has a unique invariant probability measure $\mu_\infty\in \scr P_2$ such that for some constants $t_0, c_0,C>0$ we have the log-Harnack inequality
\beq\label{LHI} P_{t_0} (\log f)(\nu)\le \log P_{t_0} f (\mu)+ c_0 \W_2(\mu,\nu)^2,\ \ \mu,\nu\in \scr P_2,\end{equation} and the Talagrand inequality
\beq\label{TLI} \W_2(\mu,\mu_\infty)^2\le C \Ent(\mu|\mu_\infty),\ \ \mu\in \scr P_2\\f\in \B_b(\R^d).\end{equation}
\beg{enumerate} \item[$(1)$] If there exist constants $c_1,\ll, t_1\ge 0$ such that
\beq\label{EW} \W_2(P_t^* \mu,\mu_\infty)^2
\le c_1\e^{-\ll t} \W_2(\mu,\mu_\infty)^2,\ \ t\ge t_1, \mu\in \scr P_2, \end{equation}
then
\beq\label{ET}\beg{split} & \max\big\{c_0^{-1} \Ent(P_t^*\mu|\mu_\infty), \W_2(P_t^*\mu, \mu_\infty)^2\big\} \\
&\le  c_1\e^{-\ll (t-t_0)} \min\big\{\W_2(\mu,\mu_\infty)^2,  C\Ent(\mu|\mu_\infty)\big\},\ \ t\ge t_0+t_1, \mu\in \scr P_2.\end{split}
 \end{equation}
\item[$(2)$]  If for some constants $\ll, c_2,t_2>0$
\beq\label{EW'} \Ent (P_t^* \mu|\mu_\infty)\le c_2\e^{-\ll t} \Ent(\mu|\mu_\infty),\ \ t\ge t_2, \nu\in \scr P_2, \end{equation}
then
\beq\label{ET'}\beg{split} &\max\big\{\Ent(P_t^*\mu,\mu_\infty), C^{-1}  \W_2(P_t^* \mu,\mu_\infty)^2\big\}\\
&\le  c_2\e^{-\ll (t-t_0)} \min\big\{c_0\W_2(\mu,\mu_\infty)^2, \Ent(\mu|\mu_\infty)\big\},\ \ t\ge t_0+t_2, \mu\in \scr P_2.\end{split} \end{equation}
\end{enumerate}
  \end{thm}

When $\si\si^*$ is invertible and does not depend on the distribution, the log-Harnack inequality \eqref{LHI} has been established in \cite{W18}.
The Talagrand inequality was first found  in \cite{TAL} for $\mu_\infty$ being the Gaussian measure, and extended in \cite{BGL} to
$\mu_\infty$ satisfying the log-Sobolev inequality
\beq\label{LS0} \mu_\infty(f^2\log f^2)\le C \mu_\infty(|\nn f|^2),\ \ f\in C_b^1(\R^d), \mu_\infty(f^2)=1,\end{equation}
see  \cite{OV} for an earlier result under a curvature condition, and  see \cite{W04} for further extensions.

To illustrate this result, we consider   the  granular media type equation for probability density functions $(\rr_t)_{t\ge 0}$ on $\R^d$:
\beq\label{E01} \pp_t \rr_t=   {\rm div} \big\{ a\nn\rr_t + \rr_t a\nn (V + W\circledast\rr_t)\big\}, \end{equation}
 where   $W\circledast \rr_t$ is in \eqref{AOO}, and the functions
 $$ a: \R^d\to \R^d\otimes \R^d,\ \ V:\R^d\to\R,\ \  W:\R^d\times\R^d\to \R$$ satisfy the following assumptions.

\beg{enumerate} \item[$(H_1)$]  $a:=(a_{ij})_{1\le i,j\le d} \in C_b^2(\R^d\to \R^d\otimes\R^d)$,  and
$a\ge \ll_aI_{d}$ for some constant $\ll_a>0$.
\item[$(H_2)$] $V\in C^2(\R^d), W\in C^2(\R^d\times\R^d)$ with  $W(x,y)=W(y,x)$, and there exist   constants $\kk_0\in\R$ and $\kk_1,\kk_2,\kk_0'>0$ such that
\beq\label{ICC1}  \Hess_V\ge \kk_0I_{d},\ \ \kk_0' I_{2d}\ge \Hess_{W}\ge \kk_0 I_{2d},\end{equation}
\beq\label{ICC2} \<x, \nn V(x)\> \ge \kk_1|x|^2-\kk_2,\ \ x\in \R^d.\end{equation}
Moreover, for any $\ll>0$,
\beq\label{ICC3} \int_{\R^d\times\R^d} \e^{-V(x)-V(y)-\ll W(x,y)} \d x\d y<\infty.\end{equation}
\item[$(H_3)$] There exists a  function $b_0\in L^1_{loc}([0,\infty))$   with
$$r_0:=  \ff {\|\Hess_W\|_\infty} 4 \int_0^\infty \e^{\ff 1 4 \int_0^t b_0(s)\d s } \d t < 1$$
such that for any $x,y,z\in\R^d$,
\beg{align*}  \big\<y-x, \nn V(x)-\nn V(y) +\nn W(\cdot,z)(x)-\nn W(\cdot, z)(y)\big\>
  \le   |x-y| b_0(|x-y|).\end{align*} \end{enumerate}

\

 For any $N\ge 2$, consider the Hamiltonian for the   system of $N$ particles:
  $$H_N(x_1,\cdots, x_N)=\sum_{i=1}^N  V(x_i)+ \ff 1 {N-1} \sum_{1\le i<j\le N}^N W(x_i, x_j),$$
  and the corresponding finite-dimensional Gibbs measure
 $$\mu^{(N)}(\d x_1,\cdots, x_N)= \ff 1 {Z_N} \e^{-H_N(x_1,\cdots, x_N)}\d x_1\cdots\d x_N,$$
 where   $Z_N:=\int_{\R^{dN}} \e^{-H_N(x)}\d x <\infty$   due to \eqref{ICC3} in $(H_2)$.
 For any  $1\le i\le N$,  the  conditional marginal of   $\mu^{(N)}$  given $z\in \R^{d(N-1)}$ is  given by
\beg{align*} &\mu^{(N)}_z(\d x) := \ff 1 {Z_N(z)}  \e^{-H_N(x|z)} \d x,  \ \ Z_N(z):= \int_{\R^d} \e^{-H_N(x|z)} \d x,\\
&H_N(x|z):= V(x)-\log \int_{\R^{d(N-1)}} \e^{-\sum_{i=1}^{N-1} \{V(z_i) +\ff 1 {N-1} W(x, z_i)\}}\d z_1\cdots \d z_{N-1}.\end{align*}
We have the following result.

\beg{thm}\label{C2.0} Assume $(H_1)$-$(H_3)$. If there is a constant $a>0$ such that the uniform   log-Sobolev inequality
\beq\label{LSS} \mu^{(N)}_z(f^2\log f^2)\le \ff 1 \bb \mu^{(N)}_z(|\nn f|^2),\ \ f\in C_b^1(\R^d), \mu^{(N)}_z(f^2)=1, N\ge 2, z\in \R^{d(N-1)}\end{equation} holds,
then there exists a unique   $\mu_\infty \in \scr P_2$ and a constant $c>0$ such that
\beq\label{WU} \W_2(\mu_t,\mu_\infty)^2 +\Ent(\mu_t|\mu_\infty)\le c \e^{-\ll_a\bb (1-r_0)^2 t} \min\big\{\W_2(\mu_0, \mu_\infty)^2 +\Ent(\mu_0|\mu_\infty)\big\},\ \ t\ge 1\end{equation}
holds for   any probability density functions  $(\rr_t)_{t\ge 0}$ solving  \eqref{E01},  where $\mu_t(\d x) := \rr_t(x)\d x, t\ge 0.$  \end{thm}

This result    allows $V$ and W to be non-convex.   For instance, let $V=V_1+V_2\in C^2(\R^d)$ such that
$\| V_1\|_\infty\land\|\nn V_1\|_\infty<\infty$, $\Hess_{V_2}\ge \ll I_{d}$ for some $\ll>0$, and   $W\in C^2(\R^d\times\R^d)$ with $\|W\|_\infty\land \|\nn W\|_\infty<\infty$. Then the uniform log-Sobolev inequality
\eqref{LSS} holds for some constant $\bb>0$. Indeed, by the Bakry-Emery criterion, $\mu_2(\d x):=\ff 1 {\int_{\R^d} \e^{-V(x)}\d x}  \e^{-V_2(x)}\d x$ satisfies the log-Sobolev inequality
$$\mu_2(f^2\log f^2)\le \ff 2 \ll \mu_2(|\nn f|^2),\ \ f\in C_b^1(\R^d),  \mu_2(f^2)=1.$$
then \eqref{LSS}  with some constant $\bb>0$ follows by the stability of the log-Sobolev inequality under bounded perturbations (see \cite{DS,CW97}) as well as Lipschitz perturbations (see \cite{Aida})  for the potential $V_2$.
Moreover, assumptions $(H_1)$-$(H_3)$  hold provided $\|\Hess_W\|_\infty$ is small enough such that $r_0<1$. So,  Theorem \ref{C2.0} applies.
See \cite{GLW} for more concrete examples satisfying $(H_1)$-$(H_3)$  and \eqref{LSS}.

\subsection{The non-degenerate case}
In this part, we make the following assumptions:
\beg{enumerate}  \item[$(A_1)$]  $b$ is continuous on $\R^d\times\scr P_2$ and there exists a constant $K>0$ such $\eqref{KK}$ holds.  
\item[$(A_2)$]  $a:=\si\si^*$ is invertible with $\ll:=\|(\si\si^*)^{-1}\|_\infty<\infty$, and there exist constants $K_2>K_1\ge 0$ such that  for any $x,y\in \R^d$ and $\mu,\nu\in \scr P_2$,
$$ \|\si(x)-\si(y)\|_{HS}^2 + 2\<b(x,\mu)-b(y,\nu), x-y\>\le K_1 \W_2(\mu,\nu)^2 - K_2  |x-y|^2.$$
Moreover, there exists a constant $K\in\R$ such that 
$$\Big\<\ff 1 2 (\nn_b a)v -\nn_{av}b,\ v\big\> \ge  K |v|^2,\ \ v\in\R^d.$$
\end{enumerate}
According to \cite[Theorem 2.1]{W18}, if $(A_1)$ holds and $b(x,\mu)$ is continuous on $\R^d\times \scr P_2$, then for any initial value $X_0\in L^2(\OO\to\R^d,\F_0,\P)$, \eqref{E1}  has a unique solution which satisfies
$$\E\Big[\sup_{t\in [0,T]} |X_t|^2 \Big]<\infty,\ \ T\in (0,\infty).$$
Let $P_t^*\mu=\L_{X_t}$ for the solution with $\L_{X_0}=\mu.$  We have the following result.

\beg{thm}\label{T1}  Assume $(A_1)$ and $(A_2)$. Then  $P_t^*$ has a unique invariant probability measure $\mu_\infty$ such that
\beq\label{EX1} \max\big\{\W_2(P_t^*\mu,\mu_\infty)^2, \Ent(P_t^*\mu|\mu_\infty)\big\}\le  \ff{c_1}{t\land 1} \e^{-(K_2-K_1)t} \W_2(\mu,\mu_\infty)^2,\ \ t>0,  \mu\in \scr P_2\end{equation} holds for some constant $c_1>0$.
If moreover $\si\in C_b^2(\R^d\to \R^d\otimes \R^m)$, then there exists a constant $c_2>0$ such that  for any $\mu\in \scr P_2, t\ge 1$,
\beq\label{EX2} \max\big\{\W_2(P_t^*\mu,\mu_\infty)^2, \Ent(P_t^*\mu|\mu_\infty)\big\} \le  c_2  \e^{-(K_2-K_1)t} \min\big\{\W_2(\mu,\mu_\infty)^2, \Ent(\mu|\mu_\infty)\big\}.    \end{equation}
\end{thm}

 \

 To illustrate this result, we consider the granular media equation \eqref{E00}, for which we take
\beq\label{SBB1} \si=\ss 2 I_{d},~~~~~~ b(x,\mu)=  -\nn \big\{V+ W\circledast\mu\big\}(x),\ \ (x,\mu)\in \R^d\times \scr P_2.\end{equation}
 The   following example is not included by   Theorem  \ref{C2.0} since the function $W$ may be non-symmetric.

\paragraph{Example 2.1 (Granular media equation).}  Consider \eqref{E0} with $V\in C^2(\R^d)$ and $W\in  C^2(\R^d\times\R^d)$ satisfying
\beq\label{CC1} \Hess_V\ge \ll I_{d},\ \ \Hess_W \ge \dd_1 I_{d},\ \ \|\Hess_W\|\le \dd_2\end{equation}
 for some constants $\ll_1, \dd_2>0$ and $\dd_1\in \R$.
 If   $\ll+\dd_1-\dd_2>0$, then there exists a unique    $\mu_\infty\in \scr P_2$ and a
 constant $c>0$ such that for any probability density functions   $(\rr_t)_{t\ge 0}$ solving \eqref{E00},   $\mu_t(\d x):=\rr_t(x)\d x$ satisfies
\beq\label{ERR0} \max\big\{\W_2(\mu_t,\mu_\infty),  \Ent(\mu_t|\mu_\infty)\big\} \le c\e^{-(\ll+\dd_1-\dd_2)t} \min\big\{\W_2(\mu_0,\mu_\infty),  \Ent(\mu_0|\mu_\infty)\big\},\ \ t\ge 1.\end{equation}

\beg{proof}  Let  $\si$ and $b$ be in \eqref{SBB1}.  Then \eqref{CC1} implies
 $(A_1)$   and
$$\<b(x,\mu)-b(y,\nu), x-y\> \le -(\ll_1+\dd_1)|x-y|^2 + \dd_2 |x-y|\W_1(\mu,\nu),$$
where we have used the formula
$$\W_1(\mu,\nu)= \sup\{\mu(f)-\nu(f):\ \|\nn f\|_\infty\le 1\}.$$ So, by taking  $\aa=\ff{\dd_2}2$ and noting that $\W_1\le \W_2$, we obtain
\beg{align*} &\<b(x,\mu)-b(y,\nu), x-y\> \le -\big(\ll+\dd_1-\aa \big)|x-y|^2 + \ff{\dd_2^2}{4\aa}  \W_1(\mu,\nu)^2\\
&\le - \Big(\ll+\dd_1-\ff {\dd_2}2\Big) |x-y|^2+ \ff{\dd_2} 2  \W_2(\mu,\nu)^2,\ \ x,y\in\R^d,\mu,\nu\in \scr P_2.\end{align*}
Therefore, if \eqref{CC1} holds for $\ll+\dd_1-\dd_2>0$, Theorem \ref{T1} implies that $P_t^*$ has a unique invariant probability measure $\mu_\infty\in\scr P_2$, such that
\eqref{ERR0} holds for $\mu_0\in \scr P_2$.
  When $\mu_0\notin \scr P_2$, we have $\W_2(\mu_0,\mu_\infty)^2=\infty$ since $\mu_\infty\in \scr P_2$.
Combining this with the Talagrand inequality
$$\W_2(\mu_0,\mu_\infty)^2\le C \Ent(\mu_0|\mu_\infty)$$
for some constant $C>0$, see the proof of Theorem \ref{T1}, we have $\Ent(\mu_0|\mu_\infty)=\infty$ for $\mu_0\notin \scr P_2$, so that \eqref{ERR0} holds for all $\mu_0\in \scr P$.\end{proof}

 \subsection{The degenerate case}
When $\R^k$ with some $k\in \mathbb N$ is considered, to emphasize the space we use $\scr P(\R^k)$ ($\scr P_2(\R^k)$) to denote the class of probability measures (with finite second moment) on $\R^k$.
 Consider the following McKean-Vlasov stochastic Hamiltonian system for   $(X_t,Y_t)\in \R^{d_1+d_2}:= \R^{d_1}\times \R^{d_2}:$
\beq\label{E21} \beg{cases} \d X_t= BY_t\d t,\\
\d Y_t=   \ss 2 \d W_t - \Big\{ B^* \nn V(\cdot,\L_{(X_t,Y_t)})(X_t) + \bb B^* (BB^*)^{-1}X_t+Y_t\Big\}\d t,\end{cases}\end{equation}
where $\bb>0$ is a constant, $B$ is a $d_1\times d_2$-matrix such that $BB^*$ is invertible, and $$V: \R^{d_1}\times \scr P_2(\R^{d_1+d_2})\to \R^{d_2}$$ is measurable.
 Let
\beg{align*}& {\psi_B}((x,y),(\bar x,\bar y)):=\ss{|x-\bar x|^2 +|B(y-\bar y)|^2},\ \ (x,y), (\bar x, \bar y) \in \R^{d_1+d_2},\\
&\W_2^{\psi_B}(\mu,\nu):=\inf_{\pi\in \C(\mu,\nu)} \bigg\{ \int_{\R^{d_1+d_2}\times\R^{d_1+d_2}} {\psi_B}^2 \d\pi\bigg\}^{\ff 1 2},\ \ \mu,\nu\in\scr P_2(\R^{d_1+d_2}).\end{align*}
We assume
\beg{enumerate} \item[{\bf (C)}]   $V(x,\mu)$ is differentiable in $x$ such that $\nn V(\cdot,\mu)(x)$ is Lipschitz continuous in $(x,\mu)\in \R^{d_1}\times \scr P_2(\R^{d_1+d_2}).$
Moreover, there exist   constants $\theta_1, \theta_2\in \R$ with
\beq\label{AH1} \theta_1+\theta_2<\bb,\end{equation}   such that  for any $(x,y), (x',y')\in \R^{d_1+d_2}$ and $\mu,\mu'\in \scr P_2(\R^{d_1+d_2})$,
\beq\label{AH2}\beg{split} & \big\< BB^* \{\nn V(\cdot,\mu)(x)-\nn V(\cdot,\mu')(x')\},  x-x'+(1+\bb)B(y-y')\big\>\\
&\ge -\theta_1{\psi_B} ((x,y), (x',y'))^2 -\theta_2 \W_2^{\psi_B} (\mu,\mu')^2.\end{split}\end{equation}
\end{enumerate}
Obviously, {\bf (C)} implies $(A_1)$ for $d=m=d_1+d_2$,  $\si= {\rm diag}\{0,\ss 2 I_{d_2}\}$, and
$$b((x,y),\mu)= \big(By, - B^* \nn V(\cdot,\mu)(x) - \bb B^* (BB^*)^{-1}x-y\big).$$ So, according to \cite{W18},  \eqref{E21} is well-posed for  any initial value in $L^2(\OO\to \R^{d_1+d_2}, \F_0,\P)$.
Let $P_t^*\mu=\L_{(X_t,Y_t)}$ for the solution with initial distribution $\mu\in \scr P_2(\R^{d_1+d_2}).$
In this case, \eqref{E21} becomes 
 $$\beg{cases}  \d X_t= BY_t\d t,\\
\d Y_t=   \ss 2 \d W_t  +Z_t(X_t,Y_t) \d t,\end{cases}$$
where $Z_t(x,y):= -  B^* \{\nn V(\cdot, P_t^*\mu)\}(x) + \bb B^* (BB^*)^{-1}x+y.$ According to \cite[Theorems 2.4 and 3.1]{W14b}, when $\Hess_V(\cdot, P_t^*\mu)$ is bounded,  
$$\rr_t(z):=\ff{(P_t^*\mu)(\d z) }{\d z} =\ff{\d (\L_{(X_t,Y_t)})(\d z)}{\d z} $$
exists and is differentiable in $z\in \R^{d_1+d_2}$. Moreover, since {\bf (C)} implies that 
  the class $$\{\pp_{y_j}, [\pp_{ y_j}, (By)_i\pp_{x_i}]: 1\le i\le d_1, 1\le j\le d_2\}$$ spans the tangent space at any point (i.e. the H\"ormander condition of rank $1$ holds), 
according to  the H\"ormander theorem,   $\rr_t\in C^\infty(\R^{d_1+d_2})$  for $t>0$ provided  $Z_t\in C^\infty(\R^{d_1+d_2})$ for $t\ge 0$. 

\beg{thm}\label{T2} Assume {\bf (C)}. Then $P_t^*$ has a unique invariant probability measure $\mu_\infty$ such that for any  $t>0$ and $\mu\in \scr P_2(\R^{d_1+d_2}),$
\beq\label{KK0} \max\big\{\W_2(P_t^*\mu,\mu_\infty)^2, \Ent(P_t^*\mu|\mu_\infty) \big\}\le \ff{c\e^{-2\kk t}}{(1\land t)^{3}}
\min\big\{\Ent(\mu|\mu_\infty), \W_2(\mu,\mu_\infty)^2\big\}\end{equation}
holds for some constant $c>0$ and
\beq\label{KK} \kk:=\ff{ 2(\bb-\theta_1-\theta_2)}{2+2\bb+\bb^2+\ss{\bb^4+4}}>0.\end{equation}
\end{thm}

\paragraph{Example 2.2 (Degenerate granular media equation).}  Let $ m\in \mathbb N$ and $W\in C^\infty(\R^m\times\R^{2m}).$ Consider  the following PDE for probability density functions $(\rr_t)_{t\ge 0}$ on $\R^{2m}$:
\beq\label{*EN2} \pp_t\rr_t(x,y)= \DD_y \rr_t(x,y) -\<\nn_x\rr_t(x,y), y\>+ \<\nn_y \rr_t(x,y), \nn_x (W\circledast\rr_t)(x) + \bb x+y\>,\end{equation}
where $\bb>0$ is a constant,   $\DD_y, \nn_x,\nn_y$ stand for the Laplacian in $y$ and the   gradient operators in $x,y$ respectively,  and
$$(W\circledast \rr_t)(x):=\int_{\R^{2m}} W(x,z)\rr_t(z) \d z,\ \ x\in \R^m.$$   If there exists a constant $\theta\in \big(0, \ff{2\bb}{1+3\ss{2+2\bb+\bb^2}}\big)$ such that
\beq\label{*EN3} |\nn W(\cdot,z)(x)-\nn W(\cdot,\bar z)(\bar x) |\le \theta \big(|x-\bar x|+|z-\bar z|\big),\ \ x,\bar x\in \R^m, z,\bar z\in \R^{2m},\end{equation}
then  there exists a unique probability measure  $\mu_\infty\in\scr P_2( \R^{2m})$   and a constant $c>0$ such that for any probability density functions $(\rr_t)_{t\ge 0}$ solving \eqref{*EN2},  $\mu_t(\d x):=\rr_t(x)\d x$ satisfies
\beq\label{ERR} \max\big\{\W_2(\mu_t,\mu_\infty)^2, \Ent(\mu_t|\mu_\infty) \big\}   \le c\e^{-\kk t} \min\big\{\W_2(\mu_0,\mu_\infty)^2, \Ent(\mu_0|\mu_\infty)\big\},\ \ t\ge 1\end{equation} holds for
$\kk= \ff{2\bb- \theta\big(1+ 3\ss{2+2\bb+\bb^2}\big)}{2+2\bb+\bb^2+\ss{\bb^4+4}}>0.$

\beg{proof} 
Let $d_1=d_2=m$ and $(X_t,Y_t)$ solve \eqref{E21} for
\beq\label{GPP} B:= I_{m},\ \ V(x, \mu):=  \int_{\R^{2m}} W(x,z)\mu(\d z).\end{equation}
We first observe that $\rr_t$ solves \eqref{*EN2} if and only if $\rr_t(z)=\ff{\d (P_t^*\mu)(\d z)}{\d z} $ for $\mu(\d z)=\rr_0(z)\d z$, where $P_t^*\mu:=\L_{(X_t, Y_t)}.$
 
Firstly,   let  $\rr_t(z) = \ff{\L_{(X_t, Y_t)}(\d z)} {\d z}$ which exists and is smooth as explained before Theorem \ref{T2}.  
By It\^o's formula and the integration by parts formula, for any $f\in C_0^2(\R^{2m})$ we have
\beg{align*}& \ff{\d}{\d t} \int_{\R^{2m}} ( \rr_t f)(z)\d z =\ff{\d }{\d t} \E [f(X_t,Y_t)]\\
&= \int_{\R^{2m}}\rr_t(x,y)\big\{\DD_y f(x,y)+\<\nn_x f(x,y), y\> -\<\nn_y f(x,y), \nn_x V(x, \rr_t(z)\d z)+\bb x+y\>\big\}\d x\d y\\
&= \int_{\R^{2m}} f(x,y) \big\{\DD_y\rr_t (x,y)-  \<\nn_x \rr_t(x,y), y\> +\<\nn_y \rr_t(x,y), \nn_x \mu_t(W(x,\cdot)) +\bb x+y\>\big\}\d x\d y.
\end{align*} Then $\rr_t$ solves \eqref{*EN2}.

On the other hand, let $\rr_t$ solve \eqref{*EN2} with $\mu_0(\d z):= \rr_0(z)\d z\in \scr P_2( \R^{2m})$. By the integration by parts formula,  $\mu_t(\d z):=\rr_t(z)\d z$ solves the non-linear Fokker-Planck equation
$$\pp_t\mu_t=L_{\mu_t}^* \mu_t$$ in the sense that for any $f\in C_0^\infty(\R^{d_1+d_2})$ we have
$$\mu_t(f)= \mu_0(f)+ \int_0^t \mu_s(L_{\mu_s}f)\d s,\ \ t\ge 0,$$
where $L_\mu:= \DD_y +y\cdot\nn_x -\{\nn_x \mu(W(x,\cdot)) + \bb x-y\}\cdot\nn_y.$ By the superposition principle, see \cite[Section 2]{BR},
we have $\mu_t=P_t^*\mu$. 

Now, as explained in the proof of Example 2.1, by Theorem \ref{T2} we only need to verify {\bf (C)} for $B, V$ in \eqref{GPP}  and
\beq\label{TTH} \theta_1= \theta\Big(\ff 1 2 +\ss{2+2\bb+\bb^2}\Big),\ \ \theta_2= \ff \theta 2 \ss{2+2\bb+\bb^2},\end{equation}
so that the desired assertion holds for
$$ \kk:=\ff{ 2(\bb-\theta_1-\theta_2)}{2+2\bb+\bb^2+\ss{\bb^4+4}}=  \ff{2\bb- \theta(1+ 3\ss{2+2\bb+\bb^2})}{2+2\bb+\bb^2+\ss{\bb^4+4}}.$$
By \eqref{*EN3} and $V(x,\mu):=\mu(W(x,\cdot))$, for any constants $\aa_1,\aa_2,\aa_3>0$ we have
\beg{align*} I &:= \big\<\nn V(\cdot,\mu)(x)- \nn V(\cdot,\bar\mu)(\bar x), x-\bar x+(1+\bb)(y-\bar y)\big\> \\
&= \int_{\R^{2m}} \big\< \nn W(\cdot, z)(x)- \nn W(\cdot, z)(\bar x), x-\bar x +(1+\bb)(y-\bar y)\big\>\mu(\d z) \\
&\qquad +\big\<\mu(\nn_{\bar x}W(\bar x,\cdot))-\bar\mu(\nn_{\bar x} W(\bar x,\cdot)),  x-\bar x +(1+\bb)(y-\bar y)\big\>\\
&\ge -\theta\big\{ |x-\bar x|+ \W_1(\mu,\bar \mu)\big\} \cdot\big(|x-\bar x|+(1+\bb)|y-\bar y|\big)\\
&\ge - \theta (\aa_2+\aa_3) \W_2(\mu,\bar\mu)^2- \theta\Big\{\Big(1+\aa_1+\ff 1 {4\aa_2}\Big)|x-\bar x|^2 +(1+\bb)^2 \Big(\ff 1 {4 \aa_1}+\ff 1 {4\aa_3}\Big)|y-\bar y|^2\Big\}.\end{align*}
Take
$$\aa_1= \ff{\ss{2+2\bb+\bb^2}-1} 2,\ \ \aa_2= \ff 1 {2\ss{2+2\bb+\bb^2 }},\ \ \aa_3=\ff {(1+\bb)^2}{2\ss{2+2\bb+\bb^2 }}.$$
We have
\beg{align*}&1+\aa_1+\ff 1 {4\aa_2} =\ff 1 2+ \ss{2+2\bb+\bb^2},\\
&(1+\bb)^2\Big(\ff 1 {4 \aa_1}+\ff 1 {4\aa_3}\Big)=\ff 1 2+ \ss{2+2\bb+\bb^2},\\
& \ \aa_2+\aa_3= \ff 1 2 \ss{2+2\bb+\bb^2}.\end{align*}
Therefore,
 $$I\ge -\ff \theta 2 \ss{2+2\bb+\bb^2} \W_2(\mu,\bar\mu)^2 -\theta\Big(\ff 1 2 + \ss{2+2\bb+\bb^2}\Big) |(x,y)-(\bar x,\bar y)|^2,$$
 i.e. {\bf (C)} holds for $B$ and $V$ in \eqref{GPP} where $B=I_m$ implies that ${\psi_B}$ is the Euclidean distance on $\R^{2m}$, and for    $\theta_1,\theta_2$ in \eqref{TTH}.
\end{proof}

 \section{Proofs of Theorems \ref{T0} and   \ref{C2.0}}

 \beg{proof}[Proof of  Theorem \ref{T0}]  (1) Since
 $$\Ent(P_{t_0}^*\nu|P_{t_0}^*\mu)=\sup_{f\ge 0, (P_{t_0}f)(\mu)=1} P_{t_0}(\log f)(\nu),$$
 \eqref{LHI} implies
 $$\Ent(P_{t_0}^*\nu|P_{t_0}^*\mu)\le c_0 \W_2(\mu,\nu)^2.$$
 This together with  $P_{t_0}^*\mu_\infty=\mu_\infty$ gives
 \beq\label{EWW} \Ent(P_{t_0}^*\mu| \mu_\infty)\le c_0 \W_2(\mu,\mu_\infty)^2,\ \ \mu\in \scr P_2.\end{equation}
 Combining   \eqref{EW} with \eqref{TLI} and \eqref{EWW},  we obtain
 $$\W_2(P_t^*\mu, \mu_\infty)^2 \le c_1\e^{-\ll t} \W_2(\mu,\mu_\infty)^2\le c_1\e^{-\ll t} \min\big\{\W_2(\mu,\mu_\infty)^2, C \Ent(\mu|\mu_\infty)\big\},\ \ t\ge t_1$$ and
  \beg{align*} &\Ent(P_t^*\mu|\mu_\infty)  =\Ent(P_{t_0}^*P_{t-t_0}^*\mu|\mu_\infty)\le c_0 \W_2(P_{t-t_0}^*\mu, \mu_\infty)^2
 \le c_0c_1\e^{-\ll (t-t_0)} \W_2(\mu, \mu_\infty)^2\\
 &= \{c_0c_1\e^{\ll t_0} \}\e^{-\ll t} \min\big\{\W_2(\mu, \mu_\infty)^2, C\Ent(\mu|\mu_\infty)\big\},\ \ t\ge t_0+t_1.\end{align*}Therefore,
 \eqref{ET} holds.

 (2) Similarly, if \eqref{EW'} holds, then   \eqref{TLI} and \eqref{EWW} imply
 \beg{align*} & \Ent(P_t^*\mu|\mu_\infty)\le c_2\e^{-\ll(t-t_0)} \min\big\{\Ent(P_{t_0}^*\mu|\mu_\infty), \Ent(\mu|\mu_\infty)\big\}\\
 & \le c_2\e^{-\ll(t-t_0)} \min\big\{ c_0 \W_2(\mu,\mu_\infty)^2, \Ent(\mu|\mu_\infty)\big\},\ \ t\ge t_0+t_2\end{align*} and
   \beg{align*} &C^{-1} \W_2(P_t^*\mu,\mu_\infty)^2\le   \Ent(P_{t-t_0}^*P_{t_0}^*\mu|\mu_\infty)\\
   &\le  c_2 \min\big\{\e^{-\ll t}  \Ent( \mu|\mu_\infty),  \e^{-\ll (t-t_0)}  \Ent( P_{t_0}^*\mu|\mu_\infty) \big\}\\
   &\le  c_2 \e^{-\ll (t-t_0)} \min\big\{ \Ent( \mu|\mu_\infty), c_0  \W_2(\mu, \mu_\infty)^2\big\},\ \ t\ge t_0+t_2.\end{align*}
 Then \eqref{ET'} holds, and the proof is finished. \end{proof}

 \beg{proof}[Proof of  Theorem  \ref{C2.0}]  By \cite[Theorem 10]{GLW}, there exits a unique   $\mu_\infty\in\scr P_2$ such that
 \beq\label{MIN} \Ent^{V,W}(\mu_\infty)=0.\end{equation}  Let $\mu_0=\rr_0\d x\in \scr P_2$. We first note that $\mu_t= P_t^*\mu_0:=\L_{X_t}$ for $X_t$ solving the distribution dependent SDE \eqref{EP} with
 \beq\label{SBB} \si(x,\mu)= \ss{2a(x)},\ \ b(x,\mu)=  \sum_{j=1}^d\pp_j a_{\cdot,j}(x) -a\nn\{V+W\circledast \mu\}(x),\ \ x\in \R^d, \mu\in \scr P_2.\end{equation} Obviously, for this choice of $(\si,b)$, assumptions $(H_1)$ and $(H_2)$ imply condition \eqref{KK} for some constant $K>0$,
 so that the McKean-Vlasov  \eqref{EP} SDE is weakly and strongly well-posed.
 For any $N\ge 2$, let $\mu_t^{(N)}=\L_{X_t^{(N)}}$ for the mean field particle system $X^{(N)}_t=(X_t^{N,k})_{1\le i\le N}$:
 \beq\label{MF} \d X_t^{N,k} = \ss 2 \si(X_t^{N,k})  \d B_t^k +\Big\{\sum_{j=1}^d  \pp_j a_{\cdot,j}(X_t^{N,k}) -a(X_t^{N,k}) \nn_k  H_N(X_t^{(N)} )\Big\} \d t,\ \  t\ge 0,\end{equation}
 where $\nn_k$ denotes the gradient in the $k$-th component, and $\{X_0^{N,k}\}_{1\le i\le N}$ are i.i.d. with distribution $\mu_0\in \scr P_2$.  According to the propagation of chaos, see \cite{SN},
  $(H_1)$-$(H_3)$ imply
 \beq\label{SNN}  \lim_{N\to \infty}\W_2(\L_{X_t^{N,1}}, P_t^*\mu_0)=0.\end{equation}

Next,  our conditions imply (25) and (26) in \cite{GLW} for $\rr_{LS}= \bb (1-r_0)^2.$ So,   by \cite[Theorem 8(2)]{GLW},
  we have the log-Sobolev inequality
 \beq\label{LSN} \mu^{(N)}(f^2\log f^2)\le \ff 2 {\bb (1-r_0)^2}\mu^{(N)}(|\nn f|^2),\ \ f\in C_b^1(\R^{dN}), \mu^{(N)}(f^2)=1.\end{equation}
 By \cite{BGL},   this implies the Talagrand inequality
 \beq\label{TLN} \W_2(\nu^{(N)}, \mu^{(N)})^2\le \ff 2 {\bb (1-r_0)^2}  \Ent(\nu^{(N)}|\mu^{(N)}),\ \ t\ge 0, N\ge 2, \nu^{(N)}\in \scr P(\R^{dN}).\end{equation}
 On the other hand,  by It\^o's formula we see that the generator of the diffusion process $X_t^{(N)} $ is
 $$L^{(N)} (x^{(N)} )= \sum_{i,j,k=1}^d \Big\{a_{ij}(x^{N,k}) \pp_{x_i^{N,k}} \pp_{x_j^{N,k}} +\pp_j a_{ij} (x^{N,k})\pp_{x_i^{N,k}} -a_{ij}(x^{N,k}) \big[\pp_{x_i^{N,k}} H_N(x^{(N)})\big]\pp_{x_i^{N,k}}\Big\},$$
 for $x^{(N)}=(x^{N,1},\cdots, x^{N,N})\in \R^{dN},$ where $x_i^{N,k}$ is the $i$-th component of $x^{N,k}\in \R^d$.  Using the integration by parts formula, we see that this operator is symmetric in $L^2(\mu^{(N)})$:
 $$\EE^{(N)}(f,g):=\int_{\R^{dN}} \<a^{(N)}\nn f,\nn g\> \d\mu^{(N)}= -\int_{\R^{dN}} (f L^{(N)}g) \d\mu^{(N)},\ \ f,g\in C_0^\infty(\R^{dN}), $$ where
  $a^{(N)}(x^{(N)}) :={\rm diag}\{a(x^{N,1}),\cdots, a(x^{N,N})\}, x^{(N)}=(x^{N,1},\cdots, x^{N,N})\in \R^{dN}.$
 So, the closure of the pre-Dirichelt form $(\EE^{(N)}, C_0^\infty(\R^{dN}))$ in $L^2(\mu^{(N)})$  is the Dieichlet form for the Markov semigroup $P_t^{(N)}$ of $X_t^{(N)}$.
  By $(H_1)$ we have  $a^{(N)} \ge \ll_a I_{dN}$, so that \eqref{LSN} implies
 $$\mu^{(N)}(f^2\log f^2)\le \ff 2 {\bb \ll_a (1-r_0)^2}\EE^{(N)}(f,f),\ \ f\in C_b^1(\R^{dN}), \mu^{(N)}(f^2)=1.$$
 It is well known that this log-Sobolev inequality implies the exponential convergence
  \beq\label{EXN} \beg{split}&\Ent(\mu_t^{(N)}| \mu^{(N)}) \le \e^{-\ll_a\bb (1-r_0)^2 t} \Ent(\mu_0^{(N)}|\mu^{(N)})\\
  &= \e^{-\ll_a\bb (1-r_0)^2 t}   \Ent(\mu^{\otimes N}|\mu^{(N)}),\ \ t\ge 0, N\ge 2,\end{split} \end{equation}
  see for instance \cite[Theorem 5.2.1]{BGL0}.
 Moreover, since $\Hess_V$ and $\Hess_W$ are bounded from below,    $(H_1)$ implies  that the Bakry-Emery curvature
 of the generator of $X_t^{(N)}$ is bounded by a constant.  Then according to \cite{W10},
  there exists a constant $K\ge 0$ such that the Markov semigroup $P_t^{(N)}$ of
 $X_t^{(N)}$ satisfies the log-Harnack inequality
\beq\label{PLM} P_t^{(N)} \log f(x)\le \log P_t^{(N)}f(y)+ \ff{K\rr^{(N)}(x,y)^2}{2(1-\e^{-2Kt})},\ \ 0<f\in \B_b(\R^{dN}), t>0, x,y\in \R^{dN},\end{equation}
 where $\rr^{(N)}$ is the intrinsic distance induced by the Dirichlet form $\EE^{(N)}$.  Since $a^{(N)}\ge \ll_a I_{dN}$, we have $\rr^{(N)}(x,y)^2\le \ll_a^{-1} |x-y|^2$.
 So, \eqref{PLM}  implies \eqref{LHI} for $P_t^{(N)}$ replacing $P_{t_0}$ and $c_0= \ff{K}{2\ll_a(1-\e^{-2Kt})}:$
$$P_t^{(N)} (\log f)(\nu) \le \log P_t^{(N)}f(\mu)+ \ff{K\,\W_2(\mu,\nu)^2}{2\ll_a(1-\e^{-2Kt})},\ \ 0<f\in \B_b(\R^{dN}), t>0, \mu,\nu \in \scr P_2( \R^{dN}).$$
 Thus, by Theorem \ref{T0}, \eqref{EXN} implies
 \beq\label{EXN2} \W_2(\mu_t^{(N)}| \mu^{(N)})^2\le \ff{c_1  \e^{-\ll_a\bb(1-r_0)^2 t}  }{1\land t}  \W_2(\mu^{\otimes N},  \mu^{(N)})^2,\ \  t>0, N\ge 2 \end{equation}
 for some constant $c_1>0$. Moreover,   \eqref{TLN},  \eqref{MIN}  and \cite[Lemma 17]{GLW}  yield
\beq\label{EXN3} \beg{split} \lim_{N\to\infty} \ff 1 N \W_2(\mu_\infty^{\otimes N}, \mu^{(N)})^2 &\le  \limsup_{N\to\infty} \ff {2 }{\bb(1-r_0)^2 N} \Ent (\mu_\infty^{\otimes N}| \mu^{(N)})^2\\
&=\ff{2}{\bb (1-r_0)^2} \Ent^{V,W}(\mu_\infty) =0.\end{split} \end{equation}
Combining this with \eqref{EXN2} we derive
\beq\label{EXN5}  \beg{split} &\limsup_{N\to\infty} \ff 1 N \W_2(\mu_t^{(N)}, \mu_\infty^{\otimes N})^2= \limsup_{N\to\infty} \ff 1 N \W_2(\mu_t^{(N)}, \mu^{(N)})^2\\
 &\le \ff{ c_1\e^{-\ll_a\bb(1-r_0)^2 t}  }{1\land t} \limsup_{N\to\infty} \ff 1 N \W_2(\mu^{\otimes N}_0, \mu^{(N)})^2\\
 &=  \ff{ c_1 \e^{-\ll_a\bb(1-r_0)^2 t}}{1\land t}   \limsup_{N\to\infty} \ff 1 N \W_2(\mu^{\otimes N}_0, \mu_\infty^{\otimes N})^2
 = \ff{ c_1 \e^{-\ll_a\bb(1-r_0)^2 t}}{1\land t}    \W_2(\mu_0, \mu_\infty)^2\ ,\ \ t>0.\end{split}\end{equation}
 Now, let $\xi=(\xi_i)_{1\le i\le N}$ and $\eta=(\eta_i)_{1\le i\le N}$ be random variables on $\R^{dN}$ such that
 $\L_\xi= \mu_t^{(N)}, \L_\eta=\mu_\infty^{\otimes N}$ and
 $$\sum_{i=1}^N \E |\xi_i-\eta_i|^2=\E|\xi-\eta|^2= \W_2(\mu_t^{(N)},\mu_\infty^{\otimes N})^2.$$
 We have $\L_{\xi_i}= \L_{X_t^{N,1}}, \L_{\eta_i}=\mu_\infty$ for any $1\le i\le N$, so that
\beq\label{EXX} N \W_2(\L_{X_t^{N,1}},\mu_\infty)^2 \le \sum_{i=1}^N \E |\xi_i-\eta_i|^2= \W_2(\mu_t^{(N)},\mu_\infty^{\otimes N})^2.\end{equation}
 Substituting this into \eqref{EXN5}, we arrive at
 $$\limsup_{N\to\infty} \W_2(\L_{X_t^{N,1}},\mu_\infty)^2 \le   \ff{ c_1 \e^{-\ll_a\bb(1-r_0)^2 t} }{1\land t}  \W_2(\mu, \mu_\infty)^2\ ,\ \ t>0. $$
This and \eqref{SNN} imply
 \beq\label{EXN6} \W_2(P_t^*\mu,\mu_\infty)^2\le  \ff{ c_1 \e^{-\ll_a\bb(1-r_0)^2 t}}{1\land t}   \W_2(\mu, \mu_\infty)^2\ ,\ \ t>0. \end{equation}
By \cite[Theorem 4.1]{W18},   $(H_1)$-$(H_3)$  imply the log-Harnack inequality
\beq\label{LHK} P_t(\log f)(\nu)\le \log P_tf(\mu)+\ff {c_2}{1\land t} \W_2(\mu,\nu)^2,\ \ \mu,\nu\in \scr P_2, t>0\end{equation}  for some constant $c_2>0$.
Similarly to the proof of \eqref{EXX} we have
$$N\W_2(\mu_\infty, \mu^{(N,1)})^2\le \W_2(\mu_\infty^{\otimes N},\mu^{(N)})^2,$$
where $\mu^{(N,1)}:=\mu^{(N)}(\cdot\times \R^{d(N-1)})$ is the first marginal distribution of $\mu^{(N)}$. This together with \eqref{EXN3} implies
$$\lim_{N\to\infty} \W_2(\mu^{(N,1)}, \mu_\infty)^2=0.$$
Therefore, applying \eqref{LSN} to $f(x)$ depending only on the firs component $x_1$,  and letting $N\to\infty$, we derive the log-Sobolev inequality
$$\mu_\infty(f^2\log f^2)\le \ff 2 {\bb (1-r_0)^2} \mu_\infty(|\nn f|^2),\ \ f\in C_b^1(\R^d), \mu_\infty(f^2)=1.$$
By \cite{BGL}, this implies \eqref{TLI} for $C=\ff 2 {\bb (1-r_0)^2}.$
Combining this with the log-Harnack inequality and  \eqref{EXN6}, by Theorem \ref{T0} we  prove  \eqref{WU} for some constant $c>0$ and $\mu_t=\L_{X_t}=P_t^*\mu_0$ for
solutions to \eqref{EP} with $b,\si$ in  \eqref{SBB}.

Similarly to the link  of \eqref{E00} and \eqref{E0'} shown in Introduction,   for any probability density functions $\rr_t$ solving \eqref{E01}, we have $\rr_t\d x= P_t^*\mu_0$ for $\mu_0=\rr_0\d x\in \scr P_2$.  So, we have proved \eqref{WU} for $\rr_t$ solving \eqref{E01} with $\mu_0\in \scr P_2$. As explained in the proof of Example 2.1 that $\Ent(\mu_0,\mu_\infty)=\W_2(\mu,\mu_\infty)=\infty$ for $\mu_0\notin \scr P_2$, so that the desired inequality \eqref{WU} tribally true.  Then the proof is finished.
\end{proof}

 \section{Proof of Theorems \ref{T1}   }

 According to \cite[Theorem 3.1]{W18}, $(A_1)$ and $(A_2)$ imply that $P_t^*$ has a unique invariant probability measure $\mu_\infty$ and
\beq\label{001}  \W_2(P_t^*\mu,  \mu_\infty) \le \e^{-\frac{1}{2}(K_2-K_1) t} \W_2(\mu,\mu_\infty),\ \ t\ge 0, \mu \in \scr P_2,\end{equation}
 while  \cite[Corollary 4.3]{W18} implies
$$ \Ent(P_t^*\mu| \mu_\infty) \le \ff {c_0}{1\land t} \W_2(\mu,\mu_\infty)^2,\ \ t>0, \mu\in \scr P_2$$  for some constant $c_0>0$.
 Then for any $p>1$,  combining these with  $P_t^*=P_{1\land t}^* P_{(t-1)^+}^*$, we obtain
 \beg{align*} &\Ent(P_t^*\mu| \mu_\infty) =\Ent(P_{1\land t}^*P_{(t-1)^+}^* \mu| \mu_\infty) \le  \ff{c_0}  {1\land t} \W_2(P_{(t-1)^+}^*\mu,\mu_\infty)^2\\
 &\le \ff{c_0\e^{-(K_2-K_1)(t-1)^+}}  {1\land t} \W_2(\mu,\mu_\infty)^2= \ff{c_0\e^{K_2-K_1}}  {1\land t} \e^{-(K_2-K_1)t} \W_2(\mu,\mu_\infty)^2.\end{align*}
 This together  with \eqref{001} implies  \eqref{EX1}   for some constant $c_1>0$.

Now, let $\si\in C_b^2(\R^d\to \R^d\otimes \R^m)$.  To deduce \eqref{EX2}   from \eqref{EX1}, it remains to find a constant $c>0$ such that the following Talagrand inequality holds:
 $$\W_2(\mu,\mu_\infty)^2\le c\, \Ent(\mu|\mu_\infty),\ \ \mu\in \scr P_2.$$
 According to \cite{BGL}, this inequality follows from the log-Sobolev inequality
 \beq\label{LSII} \mu_\infty(f^2\log f^2)\le c \mu_\infty(|\nn f|^2),\ \ f\in C_b^1(\R^d),  \mu_\infty(f^2)=1.\end{equation}
 To prove this inequality,  we consider the diffusion process $\bar X_t$ on $\R^d$ generated by
 $$\bar L:= \ff 1 2 \sum_{i,j=1}^d (\si\si^*)_{ij}\pp_i\pp_j  +\sum_{i=1}^\infty b_i(\cdot,\mu_\infty)\pp_i,$$
 which can be constructed by solving the SDE
 \beq\label{BSDE} \d\bar X_t= \si(\bar X_t)\d W_t+ b(\bar X_t,\mu_\infty)\d t.\end{equation}
 Let $\bar P_t$ be the associated Markov semigroup.  Since $P_t^*\mu_\infty=\mu_\infty$, when $\L_{X_0}=\mu_\infty$ the SDE \eqref{BSDE}   coincides with \eqref{E1}
 so that by the uniqueness, we see that $\mu_\infty$ is an invariant probability measure of $\bar P_t$. Combining this with   $(A_2)$ and It\^o's formula, we obtain
 \beq\label{ECC} \W_2(\L_{\bar X_t},\mu_\infty)^2\le \e^{-K_2t} \W_2(\L_{\bar X_0},\mu_\infty)^2,\ \ t>0. \end{equation}
  To prove the log-Sobolev inequality \eqref{LSII}, we first verify the hyperboundedness of $\bar P_t$, i.e. for   large $t>0$ we have
 \beq\label{HP} \|\bar P_t\|_{L^2(\mu_\infty)\to L^4(\mu_\infty)}<\infty.\end{equation}
 It is easy to see that conditions $(A_1)$ and $(A_2)$ in Theorem \ref{T1}  imply that $\si$ and $b(\cdot,\mu_\infty)$ satisfy conditions  $(A1)$-$(A3)$ in \cite{W11} for
 $K=-(K_2-K_1), \ll_t^2=\ll$ and $\dd_t= \|\si\|_\infty.$ So, by \cite[Theorem 1.1(3)]{W11}, we find a constant $C>0$ such that    the following  Harnack inequality  holds:
 $$(\bar P_tf(x))^2\le \bar P_t f^2 (y) \exp\Big[\ff{C|x-y|^2}{\e^{(K_2-K_1)t}-1}\Big],\ \ t>0.$$
 Then for any $f$ with $\mu_\infty(f^2)\le 1$, we have
 \beg{align*} &\big(\bar P_t f(x)\big)^2 \int_{\R^d} \exp\Big[-\ff{C|x-y|^2}{\e^{(K_2-K_1)t}-1}\Big]\mu_\infty(\d y)\\
 &\le \mu_\infty(\bar P_t f^2) = \mu_\infty(f^2)\le 1.\end{align*}
 So,
 \beq\label{UPP} \beg{split} &\sup_{\mu_\infty(f^2)\le 1} |\bar P_t f(x)|^4\le  \ff 1 {\big(\int_{\R^d} \e^{-\ff{C|x-y|^2}{\e^{(K_2-K_1)t}-1}}\mu_\infty(\d y)\big)^2} \\
&\le \ff 1 {\big(\int_{B(0,1)} \e^{-\ff{C|x-y|^2}{\e^{(K_2-K_1)t}-1}}\mu_\infty(\d y)\big)^2}\le C_1\exp\big[C_1 \e^{-(K_2-K_1)t}|x|^2\big],\ \ t\ge 1, x\in \R^d.\end{split}  \end{equation}
Next, by
 $ \|\si\|_\infty<\infty$, $(A_2)$ and It\^o's formula, for any $k\in (0,K_2)$ there exists a constant $c_k>0$ such that
$$ \d |\bar X_t|^2\le 2\<\bar X_t, \si(\bar X_t)\d W_t\> +\big\{c_k-k |\bar X_t|^2\big\}\d t.$$
Then for any $\vv>0$,
\beq\label{PW} \d\e^{\vv |\bar X_t|^2} \le  2\vv \e^{\vv|\bar X_t|^2}  \<\bar X_t, \si(\bar X_t)\d W_t\>+\vv \e^{\vv |\bar X_t|^2}\big\{c_k +2\vv\|\si\|_\infty^2|\bar X_t|^2-k|\bar X_t|^2\big\}\d t.\end{equation}
When $\vv>0$ is small enough such that $2\vv\|\si\|_\infty^2<K_2$, there exist  constants $c_1(\vv), c_2(\vv)>0$ such that
$$\vv \e^{\vv |\bar X_t|^2}\big\{c_k +2\vv\|\si\|_\infty^2|\bar X_t|^2-k|\bar X_t|^2\big\}\le c_1(\vv)-c_2(\vv)\e^{\vv|\bar X_t|^2}.$$Combining this with \eqref{PW} we obtain
$$\d \e^{\vv |\bar X_t|^2} \le c_1(\vv) - c_2(\vv) \e^{\vv |\bar X_t|^2} \d t + 2 \vv \e^{\vv |\bar X_t|^2} \<\bar X_t, \si(\bar X_t) \d W_t\>.$$
Taking for instance $\bar X_0=0$, we get
$$\ff {c_2(\vv)} t \int_0^t \E \e^{\vv|\bar X_s|^2}\d s \le \ff{1+c_1(\vv)t} t,\ \ t>0.$$
This together with \eqref{ECC} yields
$$\mu_\infty(\e^{\vv(|\cdot|^2\land N)})  = \lim_{t\to\infty} \ff 1 t\int_0^t \E \e^{\vv(|\bar X_s|^2\land N)}\d s \le \ff{c_1(\vv)}{c_2(\vv)},\ \ N>0.$$
By letting $N\to\infty$ we derive $\mu_\infty(\e^{\vv|\cdot|^2})<\infty$.
Obviously,  this and   \eqref{UPP} imply  \eqref{HP} for large $t>0$.
Moreover, by Lemma \ref{LMM} below, $(A_2)$ implies 
that for  a constant $K_0\in \R$ and all $f\in C^\infty(\R^d)$,
$$\GG_2(f):= \ff 1 2 \bar L |\si^*\nn f|^2 -\<\si^*\nn f, \si^* \nn \bar L f\> \ge K_0 |\si^*\nn f|^2,$$
i.e.  the  Bakry-Emery curvature of $\bar L$ is bounded below by  a constant $K_0$.  According to \cite[Theorem 2.1]{RW03},
this and the hyperboundedness \eqref{HP} imply the defective log-Sobolev inequality
\beq\label{DFL} \beg{split} &\mu_\infty(f^2\log f^2)\le C_1 \mu_\infty(|\si^*\nn f|^2) +C_2\\
&\le C_1\|\si\|_\infty^2 \mu_\infty(|\nn f|^2) +C_2,\ \ f\in C_b^1(\R^d), \mu_\infty(f^2)=1\end{split}\end{equation} for some constants $c_1,c_2>0$.
Since $\bar L$ is elliptic, the invariant probability measure $\mu_\infty$ is equivalent to the Lebesgue measure, see for instance \cite[Theorem 1.1(ii)]{BRW01}, so that the Dirichlet form
$$\scr E(f,g):= \mu_\infty(\<\nn f,\nn g\>),\ \ f,g\in W^{1,2}(\mu)$$ is irreducible, i.e. $f\in W^{1,2}(\mu)$ and $\scr E(f,f)=0$ imply that $f$ is  constant. Therefore,
by \cite[Corollary 1.3]{W14}, see also \cite{Miclo}, the defective log-Sobolev inequality \eqref{DFL}  implies the desired log-Sobolev inequality \eqref{LSII} for some constant $c>0$.
Hence, the proof is finished.

\beg{lem}\label{LMM} If $a:=\si\si^*\in C_b^2(\R^d;\R^{d\otimes d})$ and  is uniformly elliptic, and  there exists a constant $K\in\R$ such that 
\beq\label{ACC} \Big\<\ff 1 2 (\nn_b a)v -\nn_{av}b,\ v\big\> \ge  K |v|^2,\ \ v\in\R^d.\end{equation} Then   there exists a constant $K_0\in \R$ such that 
$$\ff 1 2 \bar L |\si^*\nn f|^2 -\<\si^*\nn f, \si^* \nn \bar L f\> \ge K_0 |\si^*\nn f|^2,\ f\in C^\infty(\R^d).$$
\end{lem}

\beg{proof} Since $a\in C_b^2$ is uniformly elliptic, there exists a constant $c_1\in\R$ such that $L_0:= \ff 1 2 {\rm tr a\nn^2}$ satisfies 
$$\ff 1 2 L_0 |\si^*\nn f|^2- \<a\nn f, \nn L_0 f\>\ge k_1|\si^*\nn f|^2, \ \ f\in C^\infty(\R^d).$$
So, it suffices to find a constant $c_2\in\R$ such that 
\beq\label{SPP} I_b(f):= \ff 1 2 \nn_b |\si^* \nn f|^2 -\<a\nn f, \nn(\nn_b f)\>\ge c_2|\si^*\nn f|^2,\ \ f\in C^\infty(\R^d).\end{equation} 
By the symmetry of $a:=\si\si^*$, we obtain
\beg{align*}& \ff 1 2 \nn_b |\si^* \nn f|^2= \big\<\nn_b (a^{\ff 1 2}\nn f), a^{\ff 1 2}\nn f\big\>= \big\<(\nn_b a^{\ff 1 2})\nn f, a^{\ff 1 2}\nn f\big\>+\Hess_f(b, a\nn f)\\
&= \Hess_f(b, a\nn f)-\ff 1 2  \big\<a^{\ff 3 2}( \nn_b a ) \nn f,  \nn f\big\>.\end{align*}
Moreover,
$$\<a\nn f, \nn(\nn_b f)\>= \Hess_f(b, a\nn f)+\big\<a\nn_{\nn f} b,\nn f\big\>.$$
So, by the boundedness of $a $, \eqref{ACC},   we obtain
$$I_b(f)=   \Big<\ff 1 2   (\nn_{b} a) \nn f -  \nn_{a\nn f} b, \nn f\big\>\ge   K|\nn f|^2\ge   K\|a\|_\infty^{-1}  \<a\nn f,\nn f\>,$$
i.e. \eqref{SPP} holds for $c_2=K\|a\|_\infty^{-1}.$ Then the proof is finished. 
 
\end{proof} 

 \section{Proof of Theorem \ref{T2} }

We first establish the log-Harnack inequality for a more general model, which extends existing results derived in  \cite{GW12, BWY15} to the distribution dependent setting.

\subsection{Log-Harnack inequality}
 Consider the following McKean-Vlasov stochastic Hamiltonian system for   $(X_t,Y_t)\in \R^{d_1}\times \R^{d_2}:$
\beq\label{E'} \beg{cases} \d X_t= \big(AX_t+BY_t)\d t,\\
\d Y_t= Z((X_t,Y_t), \L_{(X_t,Y_t)}) \d t + \si \d W_t,\end{cases}\end{equation}
where $A$ is a $d_1\times d_1$-matrix, $B$ is a $d_1\times d_2$-matrix, $\si$ is a $d_2\times d_2$-matrix,
  $W_t$ is the $d_2$-dimensional Brownian motion on a complete filtration probability space $(\OO,\{\F_t\}_{t\ge 0}, \P)$, and
  $$Z: \R^{d_1+d_2}\times \scr P_2(\R^{d_1+d_2})\to\R^{d_2}$$ is measurable. We assume
\beg{enumerate} \item[{\bf (C)}]  $\si$ is invertible, $Z$ is Lipschitz continuous, and $A,B$ satisfy the following Kalman's rank condition for some $k\ge 1$:
$${\rm Rank}[A^0B,\cdots, A^{k-1} B]=d_1, \ \ A^0:= I_{d_1}. $$ \end{enumerate}
Obviously, this assumption implies $(A_1)$, so that \eqref{E'} has a unique solution $(X_t,Y_t)$ for any initial value $(X_0,Y_0)$ with
$\mu:=\L_{(X_0,Y_0)}\in \scr P_2(\R^{d_1+d_2}).$ Let $P^*_t\mu:=\L_{(X_t,Y_t)}$   and
$$(P_tf)(\mu) :=\int_{\R^{d_1+d_2} } f\d P_t^*\mu,\ \ t\ge 0, f\in \B_b(\R^{d_1+d_2}).$$ By \cite[Theorem 3.1]{W18}, the Lipschitz continuity of $Z$ implies
\beq\label{ES0} \W_2(P_t^*\mu, P_t^*\nu) \le \e^{Kt} \W_2(\mu,\nu),\ \ t\ge 0, \mu,\nu\in \scr P_2(\R^{d_1+d_2})\end{equation} for some constant $K>0.$
We have the following result.

\beg{prp}\label{P1} Assume {\bf (C)}. Then there exists a constant $c>0$ such that
\beq\label{LH} (P_T\log f)(\nu)\le \log (P_T f)(\mu) + \ff {c\e^{cT}}{T^{4k-1}\land 1} \W_2(\mu,\nu)^2,\ \ T>0, \mu,\nu\in \scr P_2(\R^{d_1+d_2}).\end{equation}
Consequently,
\beq\label{LH'} \Ent(P_T^*\nu|P_T^*\mu) \le   \ff {c\e^{cT}}{T^{4k-1}\land 1}  \W_2(\mu,\nu)^2,\ \ T>0, \mu,\nu\in \scr P_2(\R^{d_1+d_2}).\end{equation}
 \end{prp}

\beg{proof} According to \cite[Corollary 4.3]{W18}, \eqref{LH} implies \eqref{LH'}.  Below we prove \eqref{LH} by using the coupling by change of measures summarized in \cite[Section 1.1]{W13}. By the Kalman rank condition in {\bf (C)},
$$Q_T:= \int_0^T t(T-t) \e^{(T-t)A}BB^* \e^{(T-t)A^*} \d t$$
is invertible and there exists a constant $c_1>0$ such that
\beq\label{EQ} \|Q_T^{-1}\|\le \ff{c_1\e^{c_1T}}{(T\land 1)^{2k+1}},\ \ T>0,\end{equation}
see for instance \cite[Theorem 4.2(1)]{WZ13}.

Let $(X_0,Y_0), (\bar X_0,\bar Y_0)\in L^2(\OO\to\R^{d_1+d_2},\F_0,\P)$ such that $\L_{(X_0,Y_0)}=\mu, \L_{(\bar X_0,\bar Y_0)}=\nu$ and
\beq\label{E'2} \E\big(|X_0-\bar X_0|^2 + |Y_0-\bar Y_0|^2\big) = \W_2(\mu,\nu)^2.\end{equation}
Next, let $(X_t,Y_t)$ solve \eqref{E'}. Then $\L_{(X_t,Y_t)}= P_t^*\mu.$
Consider the    the modified equation with     initial value $(\bar X_0,\bar Y_0)$:
\beq\label{E'3} \beg{cases} \d \bar X_t= \big(A\bar X_t+B\bar Y_t)\d t,\\
\d \bar Y_t= \Big\{Z((X_t,Y_t), P_t^*\mu)  +\ff {Y_0-\bar Y_0}T +\ff{\d}{\d t} \big[t(T-t) B^* \e^{(T-t)A^*}v\big]\Big\} \d t + \si \d W_t,\end{cases}\end{equation}
where
\beq\label{EV} v:= Q_T^{-1} \bigg\{\e^{TA}(X_0-\bar X_0) +\int_0^T \ff{t -T}T\e^{(T-t)A} B(\bar Y_0-Y_0) \d t\bigg\}.\end{equation}
Then
\beq\label{ES1} \beg{split} \bar Y_t-Y_t &= \bar Y_0-Y_0+\int_0^t \Big\{\ff {Y_0-\bar Y_0}T +\ff{\d}{\d r} \big[r(T-r) B^* \e^{(T-r)A^*}v\big]\Big\} \d r\\
&= \ff{T-t}T (  \bar Y_0-Y_0) + t(T-t) B^* \e^{(T-t)A^*}v,\ \ t\in [0,T].\end{split} \end{equation}
Consequently, $Y_T=\bar Y_T$, and combining with    Duhamel's formula, we obtain
\beq\label{ES2}  \beg{split} &\bar X_t-X_t= \e^{tA}(\bar X_0-X_0) +\int_0^t \e^{(t-r)A} B  \Big\{\ff{T-r}T (  \bar Y_0-Y_0) + r(T-r) B^* \e^{(T-r)A^*}v\Big\}\d r\end{split}\end{equation}
for $t\in [0,T].$ This and \eqref{EV} imply
$$\bar X_T-X_T= \e^{TA} (\bar X_0-X_0) +\int_0^T \ff{T-r}T \e^{(T-r)A} B (\bar Y_0-Y_0)\d r +Q_T v=0,$$
which together with $Y_T=\bar Y_T$ observed above yields
\beq\label{CP}  (X_T,Y_T)= (\bar X_T, \bar Y_T).\end{equation}
 On the other hand, let
 $$\xi_t=\si^{-1}\Big\{\ff 1 T (Y_0-\bar Y_0)+ \ff{\d}{\d t} \Big[t(T-t) B^* \e^{(T-t)A^*} v\Big] + Z\big((X_t,Y_t), P_t^*\mu\big)- Z\big((\bar X_t, \bar Y_t), P_t^*\nu\big)\Big\},\ \ t\in [0,T].$$
 By {\bf (C)},  \eqref{ES0}, \eqref{EQ}, \eqref{EV}, \eqref{ES1}, and \eqref{ES2}, we find  a constant $c_2>0$ such that
 \beq\label{E'4} |\xi_t|^2 \le \ff{c_2}{(T\land 1)^{4k} } \e^{c_2T}\big\{|X_0-\bar X_0|^2 + |Y_0-\bar Y_0|^2 + \W_2(\mu,\nu)^2\big\},\ \ t\in [0,T].\end{equation}
So, the Girsanov theorem implies that
$$\tt W_t:= W_t+\int_0^t\xi_sds,\ \ t\in [0,T]$$ is a $d_2$-dimensional Brownian motion under the probability measure $\Q:=R\P$, where
\beq\label{RR} R:=  \e^{-\int_0^T \<\xi_t,\d W_t\> -\ff 1 2 \int_0^T |\xi_t|^2\d t}.\end{equation}
Reformulating \eqref{E'3} as
 $$\beg{cases}  \d \bar X_t= \big(A\bar X_t+B\bar Y_t)\d t,\\
\d \bar Y_t= Z((\bar X_t,\bar Y_t), P_t^*\nu)   \d t + \si \d \tt W_t, \ \ t\in [0,T],\end{cases}$$
by the weak uniqueness of  \eqref{E'} and that the distribution of $(\bar X_0,\bar Y_0)$ under $\Q$ coincides with   $\L_{(\bar X_0, \bar Y_0)}=\nu$, we obtain
$\L_{(\bar X_t, \bar Y_t)|\Q}=P_t^*\nu$ for $t\in [0,T].$ Combining this with \eqref{CP} and using the Young inequality, for any $f\in \B_b^+(\R^{d_1+d_2})$ we have
\beq\label{LH3} \beg{split} & (P_T\log f)(\nu)= \E[R \log f(\bar X_T, \bar Y_T)] =  \E[R \log f(X_T,  Y_T)] \\
&\le \log \E[f(X_T,Y_T)] +\E [R\log R]= \log (P_Tf)(\mu) + \E_\Q [\log R].\end{split}\end{equation}
By  \eqref{E'4},  and \eqref{RR},    $\tt W_t$ is a Brownian motion under $\Q$, and noting that   $\Q|_{\F_0}=\P|_{\F_0}$ and  \eqref{E'2}  imply
$$ \E_\Q\big(|X_0-\bar X_0|^2 + |Y_0-\bar Y_0|^2\big) = \W_2(\mu,\nu)^2, $$
we find a constant $c>0$ such that $$  \E_\Q [\log R] =\ff 1 2\E_\Q\int_0^T|\xi_t|^2\d t \le \ff{c\e^{cT}}{(T\land 1)^{4k-1} } \W_2(\mu,\nu)^2.$$
Therefore, \eqref{LH} follows from \eqref{LH3}.
\end{proof}
 \subsection{Proof of Theorem \ref{T2} }

 We first prove the exponential convergence of $P_t^*$ in $\W_2$.

 \beg{lem}\label{LN1} Assume {\bf (C)}. Then there exists a constant $c_1>0$ such that
 \beq\label{AC0} \W_2(P_t^*\mu,P_t^*\nu)^2\le c_1\e^{-\kk t} \W_2(\mu,\nu)^2,\ \  t\ge 0, \mu,\nu\in \scr P_2(\R^{d_1+d_2}).\end{equation}
 Consequently, $P_t^*$ has a unique invariant probability measure $\mu_\infty\in\scr P_2(\R^{d_1+d_2})$. \end{lem}

 \beg{proof}  As shown in the proof of \cite[Theorem 3.1(2)]{W18} that the second assertion follows from the first. So, it suffices to prove \eqref{AC0}.
For  \beq\label{AC1} a:=\Big(\ff{1+\bb+\bb^2}{1+\bb}\Big)^{\ff 1 2} ,\ \ r:=  a -\ff{\bb}{ a} =\ff 1 {\ss{(1+\bb)(1+\bb+\bb^2)}}\in (0,1),\end{equation}
 we define the distance
\beq\label{AC2}  \bar \psi_{B}((x,y),(\bar x,\bar y)):=\ss{a^2|x-\bar x|^2 +|B(y-\bar y)|^2+ 2 r a \<x-\bar x, B(y-\bar y)\>}
 \end{equation}  for $ (x,y), (\bar x, \bar y) \in \R^{d_1+d_2}. $ Then there exists a constant $C>1$ such that
 \beq\label{ACC} C^{-1}|(x-\bar x, y-\bar y)|\le  \bar\psi_{B}((x,y),(\bar x,\bar y))\le C|(x-\bar x, y-\bar y)|.\end{equation}
Moreover, we claim that
\beq\label{AC3}   \bar\psi_{B}((x,y),(\bar x,\bar y))^2\le \ff{2+2\bb+\bb^2 +\ss{\bb^4+4}}{2(1+\bb)}   {\psi_B}((x,y),(\bar x,\bar y))^2.\end{equation}
 Indeed, by \eqref{AC1} and \eqref{AC2},  for any $\vv>0$ we have
\beq\label{AC}   \bar\psi_{B}((x,y),(\bar x,\bar y))^2\le a^2(1+\vv) |x-\bar x|^2 + \Big(1+ \ff 1 {\vv(1+\bb)(1+\bb+\bb^2)}\Big)|B(y-\bar y)|^2.\end{equation}
Obviously,  by \eqref{AC1},
$$\vv:= \ff{1-a^2+\ss{(a^2-1)^2+4 a^2(1+\bb)^{-1}(1+\bb+\bb^2)^{-1}}}{2 a^2}= \ff{\ss{\bb^4+4}-\bb^2}{2(1+\bb+\bb^2)}$$
satisfies
$$a^2(1+\vv)= 1+ \ff 1 {\vv(1+\bb)(1+\bb+\bb^2)}=  \ff{2+2\bb+\bb^2 +\ss{\bb^4+4}}{2(1+\bb)}.$$
Thus, \eqref{AC3} follows from \eqref{AC}.

Now, let $(X_t, Y_t)$ and $ (\bar X_t, \bar Y_t)$ solve \eqref{E21} with $\L_{(X_0,Y_0)}=\mu, \L_{(\bar X_0,\bar Y_0)}=\nu$ such that
\beq\label{AC4} \W_2(\mu,\nu)^2 =\E|(X_0-\bar X_0, Y_0-\bar Y_0)|^2.\end{equation}
Simply denote $\mu_t=\L_{(X_t,Y_t)}, \bar \mu_t= \L_{(\bar X_t, \bar Y_t)}.$ By  {\bf (C)}  and It\^o's formula, and noting that \eqref{AC1} implies
$$a^2-\bb-ra=0,\ \ 1-ra=ra\bb =\ff\bb {1+\bb},$$  we obtain
\beg{align*} &\ff 1 2 \d \big\{\bar\psi_{B}((X_t,Y_t), (\bar X_t,\bar Y_t))^2\big\}= \big\<a^2(X_t-\bar X_t)+ r a B(Y_t-\bar Y_t), B(Y_t-\bar Y_t)\big\>\d t  \\
&  + \big\<B^*B (Y_t-\bar Y_t)+ r a B^*(X_t-\bar X_t),\ \bb B^*(BB^*)^{-1} (\bar X_t-X_t) +\bar Y_t-Y_t \big\>\d t \\
&+\big\<B^*B (Y_t-\bar Y_t)+ r a B^*(X_t-\bar X_t),\ B^*\{\nn V(\bar X_t, \bar\mu_t) -\nn V(X_t,\mu_t)\}\big\> \d t\\
&\le \Big\{-(1-ra)|B(Y_t-\bar Y_t)|^2+(a^2-\bb-ra) \<X_t-\bar X_t, B(Y_t-\bar Y_t)\> - ra\bb    |X_t-\bar X_t|^2 \\
&\qquad + \big\<B^*B (Y_t-\bar Y_t)+ (1+\bb)^{-1}B^*(X_t-\bar X_t),\ B^*\{\nn V(\bar X_t, \bar\mu_t) -\nn V(X_t,\mu_t)\} \big\>\Big\}\d t\\
&\le \Big\{\ff{\theta_2}{1+\bb}\W_2^{\psi_B}(\mu_t,\bar\mu_t)^2-\ff{\bb-\theta_1}{1+\bb}{\psi_B}((X_t,Y_t), (\bar X_t,\bar Y_t))^2 \Big\}\d t.\end{align*}
By \eqref{AC3} and the fact   that
$$\W_2^{\psi_B}(\mu_t,\bar\mu_t)^2\le \E[{\psi_B}((X_t,Y_t), (\bar X_t,\bar Y_t))^2],$$
for $\kk>0$ in \eqref{KK}, we obtain
\beg{align*}&\ff 1 2 \big\{\E[ \bar\psi_{B}((X_t,Y_t), (\bar X_t,\bar Y_t))^2]- \E[ \bar\psi_{B}((X_s,Y_s), (\bar X_s,\bar Y_s))^2]\big\}\\
&\le -\ff{\bb-\theta_1-\theta_2}{1+\bb} \int_s^t \E [{\psi_B}((X_r,Y_r), (\bar X_r,\bar Y_r))^2]\d r \\
&\le - \kk  \int_s^t \E [ \bar\psi_{B}((X_r,Y_r), (\bar X_r,\bar Y_r))^2]\d r,\ \ t\ge s\ge 0.\end{align*} Therefore,  Gronwall's inequality implies
$$\E[ \bar\psi_{B}((X_t,Y_t), (\bar X_t,\bar Y_t))^2]\le \e^{-2\kk t} \E[ \bar\psi_{B}((X_0,Y_0), (\bar X_0,\bar Y_0))^2],\ \ t\ge 0.$$
Combining this with \eqref{ACC} and \eqref{AC4}, we prove \eqref{AC0} for some constant $c>0$.
\end{proof}

\beg{proof}[Proof of Theorem \ref{T2}] By Proposition \ref{P1} with $k=1$,   Lemma \ref{LN1} and Theorem \ref{T0}, we only need to verify the Talagrand inequality.
As shown in the beginning  of \cite[Section 3]{GW19} that $ \mu_\infty$ has the representation
$$\mu_\infty(\d x,\d y)= Z^{-1} \e^{\bar V(x,y)}\d x\d y,\ \ \bar V(x,y):= V(x,\mu_\infty)+ \ff\bb 2 |(BB^*)^{-\ff 1 2} x|^2 +\ff 1 2 |y|^2,$$
where $Z:=\int_{\R^{d_1+d_2} }\e^{ -\bar V(x,y)}\d x\d y$ is the normalization constant.  Since \eqref{AH2} implies
$$BB^* \Hess_{ V(\cdot,\mu_\infty)}\ge -\theta_1I_{d_1},$$
we deduce from \eqref{AH1} that
$$\Hess_{\bar V}\ge \gg I_{d_1+d_2},\ \ \gg:=1\land \ff{\bb-\theta_1}{\|B\|^2} >0.$$
So, by the Bakry-Emery criterion \cite{BE84}, we have the log-Sobolev inequality
$$\mu_\infty(f^2\log f^2)\le \ff 2 {\gg} \mu_\infty(|\nn f|^2),\ \ f\in C_b^1(\R^{d_1+d_2}), \mu_\infty(f^2)=1.$$
According to  \cite{BGL}, this implies the Talagrand inequality
$$\W_2(\mu,\mu_\infty)^2\le \ff 2 {\gg} \Ent(\mu|\mu_\infty).$$
Then the proof if finished.
\end{proof}

 \paragraph{Acknowledgement.} We would like to thank the referee for helpful comments and corrections.

\end{document}